\documentclass[12pt,a4paper]{article}

\usepackage{amsmath,amsfonts}
\usepackage{mathrsfs}
\usepackage[usenames]{color}

\newtheorem{theorem}{Theorem}

\newtheorem{proposition}{Proposition}

\def \Limsup{\mathop{\overline{\lim}}\limits}

\def\Pb{\mathbf{P}}

\def\Ex{\mathbf{E}}

\def\sgn{{\rm sgn}}
\def\1{\mbox{1\hspace{-.25em}I}}

\begin{document}

\title{On Asymptotic Distribution of Parameter Free  tests for
  Ergodic Diffusion Processes}
\author{Yury A. \textsc{Kutoyants}\\
{\small Laboratoire de Statistique et Processus, Universit\'e du Maine}\\
{\small 72085 Le Mans, C\'edex 9, France}}
\date{}

\maketitle
\begin{abstract}
We consider two problems of  constructing of  goodness of fit tests for ergodic
diffusion processes. The first one is concerned with a composite
basic hypothesis for a parametric class of diffusion processes, which  includes the
Ornstein-Uhlenbeck and simple switching processes. In this case we propose
asymptotically parameter free tests of Cram\'er-von Mises type.
The basic hypothesis in the
second problem is simple and we propose asymptotically distribution free tests
for a wider class of trend coefficients.
\end{abstract}


\section{Introduction}

In this paper we consider two different goodness of fit (GoF) hypotheses testing problems for the
diffusion process
\begin{equation*}
\label{1}
{\rm d}X_t= S\left(X_t\right)\;{\rm d}t +\sigma\left(X_t\right)\;{\rm d}W_t, \quad X_0, \quad
0\leq t\leq T.
\end{equation*}
In the first problem the observed
process under the basic hypothesis (${\cal H}_0$) satisfies the stochastic differential equation
\begin{equation*}
\label{2}
 {\rm d}X_t=-\beta \, \sgn\left(X_t-\alpha
\right)\left|X_t-\alpha\right|^\gamma {\rm   d}t+\sigma \,{\rm d}W_t,\quad
X_0,\quad 0\leq t\leq T ,
\end{equation*}
where $\vartheta =\left(\alpha ,\beta \right)\in \Theta $ is the unknown
parameter, $\beta >0, \gamma\geq 0$ and $\sigma >0$.
Therefore the hypothesis is parametric composite.

In the second problem we assume that under the basic hypothesis (${\cal H}_0$) the observed process satisfies
\begin{equation*}
\label{3}
{\rm d}X_t=S_0\left(X_t\right){\rm
  d}t+\sigma\left(X_t\right) \,{\rm d}W_t,\quad X_0,\quad 0\leq t\leq T ,
\end{equation*}
where $S_0\left(\cdot \right)$ is a known function, i.e., (${\cal H}_0$) is simple.

In both models the alternatives   are nonparametric
and, under the hypothesis ${\cal H}_0$, the diffusion processes are assumed to be ergodic with the
invariant densities $f\left(\vartheta,x\right)$ and $f_{S_0}\left(x\right)$
respectively.   We denote  the corresponding distribution functions by $F\left(\vartheta ,x\right)$ and $F_{S_0}\left(x\right)$.

Our goal is to construct the goodness of fit tests
which provide the fixed limit error $\varepsilon\in \left(0,1\right) $.
Introduce the class ${\cal K}_\varepsilon $ of such tests, i.e., the tests
$\bar\psi _T$ satisfying the relations
$$
\lim_{T\rightarrow \infty }\Ex_{\vartheta } \bar\psi _T=\varepsilon \qquad
    {\rm for \quad all}\qquad \vartheta \in \Theta ,
$$
and
$$
\lim_{T\rightarrow \infty }\Ex_{S_0} \bar\psi _T=\varepsilon
$$
in the first and the second problems respectively.

All tests studied in the present work are of the form $\hat \psi
_T=\1_{\left\{\Delta _T>c_\varepsilon \right\}}$, where $\Delta _T$ is the
Cram\'er-von Mises type statistic.  More precisely, in the first problem
$\Delta _T$ is either  of the $L_2$ distances $D\left(\hat F_T\left(x\right),
F\left(\hat\vartheta _T,x\right)\right)$ and $D\left(\hat f_T\left(x\right),
f\left(\hat\vartheta _T, x\right)\right)$, where $\hat F_T\left(x\right)$ is
the empirical distribution function,  $\hat f_T\left(x\right)$ is the local
time estimator of the invariant density  and $\hat\vartheta _T$ is the maximum
likelihood estimator (MLE) of the parameter $\vartheta $.
Similarly, in the second problem $\Delta _T$ is one of the distances $D\left(\hat F_T\left(x\right),
F_{S_0}\left(x\right)\right)$ and $D\left(\hat f_T\left(x\right),
f_{S_0}\left( x\right)\right)$.

Let us denote by $\Delta \left(\vartheta \right)$ and $\Delta
\left(S_0\right)$ the limits (in distribution) of the test statistics in the
first and the second problems. Then the thresholds $c_\varepsilon
$ in these tests have to satisfy the equations
\begin{equation}
\label{4}
\Pb_\vartheta \left\{\Delta \left(\vartheta \right)>c_\varepsilon
\right\}=\varepsilon ,\qquad \Pb_{S_0}\left\{\Delta \left(S_0
\right)>c_\varepsilon \right\}=\varepsilon  .
\end{equation}
The main contribution of this work is the following. We introduce
modifications of the statistics $\Delta _T$, so that their limit distributions  do not depend on $\vartheta $ in the
first problem and do not depend on $S_0\left(\cdot \right)$ in the second
problem. Therefore the corresponding tests are {\it asymptotically parameter
  free} in the first case and {\it asymptotically distribution free } in the
second case. These modifications essentially simplify the solution of the
equations \eqref{4} and allow  to choose the thresholds $c_\varepsilon $
before actually conducting the experiments.

Let us briefly recall  what happens in the analogous problems in the case of independent
identically distributed observations   $X_1,\ldots,X_n$. In the problem of the first type we have the
following results.
Suppose that under the basic  hypothesis
$$
{\cal H}_0\quad :\qquad \quad X_j\sim
F_0\left(\vartheta ,x\right),\quad \vartheta \in \Theta  ,
$$
where $F_0\left(\cdot ,x\right)$ is some known distribution function. The
limit distribution of the Cram\'er-von Mises staistics (under hypothesis
 ${\cal H}_0$)
$$
\Delta _n\left(X^n\right)=n\int_{-\infty }^{\infty }\left[\hat
  F_n\left({x}\right)-F_0\left(\hat\vartheta _n,x\right)\right]^2{\rm d}F_0\left(\hat\vartheta
  _n,x\right)\Longrightarrow \Delta \left(\vartheta \right)
$$
depends on $\vartheta $. Here $\hat F_n\left(x\right)$ is the empirical
distribution function and $\hat\vartheta _n$ is some estimator.  The choice of
the threshold $c_\varepsilon $ for the GoF test
$$
\hat\psi _n=\1_{\left\{\Delta _n\left(X^n\right)>c_\varepsilon \right\}}
$$
can be a difficult problem, since $c_\varepsilon =c_\varepsilon
\left(\vartheta \right)$ is solution of the equation
$$
\Pb_\vartheta \left\{\Delta \left(\vartheta \right)>c_\varepsilon
\right\}=\varepsilon .
$$

It is well-known that  for some distributions, say with shift and scale
parameters like $F\left(\frac{x-\alpha }{\beta }\right)$, this limit can be {\it
  asymptotically parameter free} (APF). For example, if the hypothesis is
$$
{\cal H}_0\qquad \quad :\quad X_j\sim
{\cal N}\left(\alpha ,\beta ^2\right),\quad \vartheta =\left(\alpha ,\beta
\right)\in \Theta ,
$$
 then the  limit distribution of the Cram\'er-von Mises staistics
$$
\Delta _n\left(X^n\right)=n\int_{-\infty }^{\infty }\left[\hat
  F_n\left({x}\right)-F_0\left(\frac{x-\hat\alpha
  _n}{\hat\beta _n}\right)\right]^2{\rm d}F_0\left(\frac{x-\hat\alpha
  _n}{\hat\beta _n}\right)
$$
does not depend on $\vartheta $ (see, e.g., \cite{KKW},\cite{Dur}, \cite{LR}).

Here
$F_0\left(x\right)$ is the  distribution function of  $N\left(0,1\right)$ random variable.
Therefore  the threshold $c_\varepsilon $ does not depend on $\vartheta $ and
the test can be easily constructed. The  similar statement
for Pareto distribution  was studied by Choulakian and Stephens \cite{CS} and
another  class of
distributions was treated by Martynov \cite{M10}.

The general case of ergodic diffusion processes with one-dimensional shift
parameter was studied by Negri and Zhou \cite{NZ}. They showed that the limit
distribution of the Cram\'er-von Mises statistic does not depend on the unknown
 parameter.

\section{Preliminaries}

We need some properties of the estimators $\hat F_T\left(x\right)$ and $\hat
f_T\left(x\right)$, which we recall below.
We assume  that the trend $S\left(x\right)$ and the diffusion $\sigma
\left(x\right)^2$ coefficients of the observed diffusion process
$$
{\rm d}X_t=S\left(X_t\right)\,{\rm d}t+\sigma \left(X_t\right)\,{\rm
  d}W_t,\quad                     X_0,                       \quad 0  \leq t\leq T
$$
satisfy the following conditions.

${\cal ES}.$ {\it The function $S\left(\cdot \right)$ is locally bounded, the
function $\sigma \left(\cdot \right)^2>0 $ is continuous  and for some $C>0$
the condition
$$
x\,S\left(x\right)+\sigma \left(x\right)^2\leq C\left(1+x^2\right)
$$
holds.}

Under this condition the stochastic differential equation has a unique weak
solution (see, e.g., \cite{Durett}).

Let us denote
$$
V\left(x\right)=\int_{0}^{x}\exp\left\{-2\int_{0}^{y}\frac{S\left(z\right)}{\sigma
  \left(z\right)^2}{\rm d}z\right\} {\rm d}y
 $$
and
$$
G\left(S \right)=\int_{-\infty }^{\infty }\sigma
\left(y\right)^{-2}\exp\left\{2\int_{0}^{x}\frac{S\left(y\right)}{\sigma
  \left(y\right)^2}{\rm d}y\right\} {\rm d}x.
$$
The next condition is:

${\cal RP}.$ {\it The functions $S\left(\cdot \right)$  and $\sigma \left(\cdot
\right)^2$ are such that }
$$
G\left(S\right)<\infty ,\qquad V\left(x\right)\longrightarrow \pm \infty \quad
{as} \qquad x\longrightarrow \pm \infty .
$$
Under this condition the diffusion process is ergodic, i.e., positive recurrent  with
the invariant density
$$
f\left(x\right)=\frac{1}{G\left(S\right)\;\sigma
  \left(x\right)^2}\;\exp\left\{2\int_{0}^{x}\frac{S\left(y\right)}{\sigma
  \left(y\right)^2}\; {\rm d}y\right\}.
$$

The empirical distribution function  $\hat F_T\left(x\right)$ and  the local
time density estimator $\hat f_T\left(x\right)$  of the invariant density
$f\left(x\right)$ are
$$
\hat F_T\left(x\right)=\frac{1}{T}\int_{0}^{T} \1_{\left\{X_t<x\right\}}\;{\rm
    d}t, \qquad \hat f_T\left(x\right)=\frac{\Lambda _T\left(x\right)}{T\sigma \left(x\right)^2}
$$
where the {\it local time} $\Lambda _T\left(x\right)$ satisfies the equation
(Tanaka-Meyer formula)
$$
\Lambda
_T\left(x\right)=\left|X_t-x\right|-\left|X_0-x\right|-\int_{0}^{T}\sgn\left(X_t-x\right)\;{\rm
  d}X_t.
$$

 Recall that these estimators  are consistent,
 asymptotically normal and asymptotically efficient  under the basic hypothesis (see \cite{Kut04}). The proof
 of these properties  is based on the representations
\begin{align}
\label{df}
\sqrt{T}\left(\hat F_T\left(x\right)-F\left( x\right)\right)&=\frac{2}{\sqrt{T}}
\int_{0}^{T}\frac{F\left(x\right)F\left(X_t\right) -F\left(x\wedge X_t\right)}{\sigma\left(X_t\right)
  \,f\left(X_t\right)}\,{\rm d}W_t \nonumber \\
&\quad +\frac{2}{\sqrt{T}}
\int_{X_0}^{X_T}\frac{F\left(y\right)F\left(x\right) -F\left(y\wedge x\right)}{\sigma\left(y\right)^2
  \,f\left(y\right)}\,{\rm d}y
\end{align}
and
\begin{align}
\label{lte}
\sqrt{T}\left(\hat f_T\left(x\right)-f\left(x\right)\right)&=\frac{2f\left(x\right)}{\sqrt{T}}
\int_{0}^{T}\frac{F\left(X_t\right)-\1_{\left\{X_t>x\right\}}}{\sigma\left(X_ty\right)
  \,f\left(X_t\right)}\,{\rm d}W_t\nonumber\\
&\quad +\frac{2f\left(x\right)}{\sqrt{T}}
\int_{X_0}^{X_T}\frac{\1_{\left\{y>x\right\}} -F\left(y\right)}{\sigma\left(y\right)^2
  \,f\left(y\right)}\,{\rm d}y.
\end{align}
Using these representations and the central limit theorem for stochastic  integrals we obtain  the limits in distribution
\begin{align*}
\sqrt{T}\left(\hat F_T\left(x\right)-F\left(x\right)\right)&\Longrightarrow
  2\int_{-\infty
  }^{\infty }\frac{F\left(y\right)F\left(x\right)-F\left(y\wedge x\right)}{\sigma
    \left(y\right)\sqrt{f\left(y\right)}}\;{\rm d}W\left(y\right),\\
\sqrt{T}\left(\hat f_T\left(x\right)-f\left(x\right)\right)&\Longrightarrow
 2f\left(x\right)\int_{-\infty }^{\infty
  }\frac{F\left(y\right)-\1_{\left\{y>x\right\}}}{\sigma
    \left(y\right)\sqrt{f\left(y\right)}}\;{\rm d}W\left(y\right),
\end{align*}
where  $W\left(\cdot \right)$ is a two-sided Wiener process.

The estimator $\hat f_T\left(x\right)$ is the  a.s.  derivative of $\hat
F_T\left(x\right)$. Indeed, using the equality (see \cite{RY})
$$
\int_{0}^{T}h\left(X_t\right)\,{\rm d}t=\int_{-\infty }^{\infty
}h\left(y\right)\,\frac{\Lambda _T\left(y\right)}{\sigma \left(y\right)^2}\; {\rm d}y
$$
we can write
\begin{align*}
\hat F_T\left(x\right)=\int_{-\infty }^{\infty
}\1_{\left\{y<x\right\}}\,\frac{\Lambda _T\left(y\right)}{T\sigma
  \left(y\right)^2}\; {\rm d}y=\int_{-\infty }^{x
}\hat f_T\left(y\right)\; {\rm d}y.
\end{align*}
As the local time is continuous with probability one  we have the limit
$$
\lim_{\alpha \rightarrow 0}\frac{\hat F_T\left(x+\alpha \right)-\hat
  F_T\left(x\right) }{\alpha } =\lim_{\alpha \rightarrow 0}\frac{1}{\alpha }\int_{x }^{x+\alpha
}\hat f_T\left(y\right)\; {\rm d}y =\hat f_T\left(x\right).
$$
Therefore we can call the  local time estimator the {\it empirical density}.
It is easy to see that the representation \eqref{lte} can be obtained from
\eqref{df} through differentiating.

Introduce the class ${\cal P}$ of locally bounded functions with polynomial
majorants ($p>0$)
$$
{\cal P}=\left\{h\left(\cdot \right):\quad \left|h\left(y\right)\right|\leq
C\left(1+\left|y\right|^p\right)\right\} .
$$
and  the following condition:

${\cal A}_0.$ {\it The functions  $S\left(\cdot \right), \sigma \left(\cdot \right)^{\pm 1} \in {\cal
    P} $ and}
$$
\Limsup_{\left|y\right|\rightarrow \infty
}\;\sgn\left(y\right)\;\frac{S\left(y\right)}{\sigma \left(y\right)^2} <0.
$$
Note that if $S\left(\cdot \right)$ and $\sigma \left(\cdot \right)$  satisfy
${\cal A}_0$ then the condition  ${\cal RP} $ is fulfilled.

Moreover, under condition ${\cal A}_0$  for any $p>0$ there exist $\kappa
>0$ and $C>0$ such that
\begin{equation*}
\Ex \left|\sqrt{T}\left(\hat f_T\left(x\right)-f\left(x\right)\right)
\right|^p\leq C\,e^{-\kappa  \left|x\right|}.
\end{equation*}
For the proof see Proposition 1.11, \cite{Kut04}.

\section{Asymptotically Parameter Free Tests}

The first problem is the following.
 We observe an ergodic diffusion process $X^T=\left(X_t,0\leq t\leq T\right)$, which
solves the  equation
\begin{equation}
\label{7}
{\rm d}X_t=S\left(X_t \right) \;{\rm d}t +\sigma
\;{\rm d}W_t,\qquad X_0, \quad 0\leq t\leq T
\end{equation}
and we have to test the composite basic hypothesis:\\

  ${\cal H}_0$ {\it this
process admits the stochastic differential }
\begin{equation}
\label{K}
{\rm d}X_t=-\beta \; \sgn\left(X_t-\alpha \right)\left|X_t-\alpha
\right|^\gamma \;{\rm d}t +\sigma
\;{\rm d}W_t,\qquad X_0, \quad 0\leq t\leq T ,
\end{equation}
{\it where $\vartheta =\left(\alpha ,\beta \right)$ is the  unknown parameter,
$\vartheta
\in \Theta =\left(a_1,a_2\right)\times \left(b_1,b_2\right), b_1>0$.}

against the nonparametric alternative:\\

 ${\cal H}_1$ {\it   the observed process
does not belong to this parametric family.}

 The parameters $\gamma\geq 0 $  and $\sigma >0$ are assumed to be  known.

Note that if $\gamma =1$, then we obtain {\it Ornstein-Uhlenbeck} process
$$
{\rm d}X_t=-\beta \; \left(X_t-\alpha \right) \;{\rm d}t +\sigma
\;{\rm d}W_t,\qquad X_0, \quad 0\leq t\leq T
$$
and if $\gamma =0$ then the solution of \eqref{K} is the {\it simple switching
  process}
$$
{\rm d}X_t=-\beta \; \sgn\left(X_t-\alpha \right) \;{\rm d}t +\sigma
\;{\rm d}W_t,\qquad X_0, \quad 0\leq t\leq T
$$
studied in \cite{Kut04}, Section 3.4. For $\gamma =3$ we have the cubic trend
$$
{\rm d}X_t=-\beta \; \left(X_t-\alpha \right)^3 \;{\rm d}t +\sigma
\;{\rm d}W_t,\qquad X_0, \quad 0\leq t\leq T .
$$

It is easy to verify that for $\beta >0, \gamma \geq 0$ this process is positive recurrent
with the invariant density
\begin{equation}
\label{pdf}
f\left(\vartheta ,x\right)=\frac{\beta ^{\frac{1}{\gamma +1}}}{G_\gamma \;\sigma
  ^{\frac{2}{\gamma +1}}}\;\exp\left\{-\frac{2\beta \left|x-\alpha
  \right|^{\gamma +1}}{\left(\gamma +1\right)\sigma ^2}\right\}.
\end{equation}
The normalizing constant is
$$
G_\gamma =\left(\frac{2}{\gamma +1}\right)^{\frac{\gamma }{\gamma+1}}\;\Gamma \left(\frac{1}{\gamma +1}\right)  ,
$$ where $\Gamma \left(\cdot \right)$ is the Gamma function. Below we denote by
$f_0\left(x\right)=f\left(\vartheta _0,x\right)$ and
$F_0\left(x\right)=F\left(\vartheta _0,x\right)$ the density and the
distribution function corresponding to the values $\vartheta
_0=\left(0,1\right), \sigma =1$ and we denote by $\xi $ the random variable with such
distribution function.

It will be convenient to study the cases $\gamma \geq {1}$ (including O-U
process) and $0\leq \gamma <\frac{1}{2}$  separately because the rates  of
convergence of
the MLE $\hat \alpha _T$ in these two cases are essentially different. We defer the discussion of
the  complementary case  $\gamma \in [\frac{1}{2},1)$ to  section 3.3 below.

\subsection{Case $\gamma \geq 1$.}

Let us consider the following ergodic diffusion process as the basic model
(under hypothesis ${\cal H}_0$)
\begin{equation*}
{\rm d}X_t=-\beta \; \sgn\left(X_t-\alpha \right) \;\left|X_t-\alpha
\right|^{\gamma }{\rm d}t+\sigma \;{\rm
  d}W_t,\quad X_0,\quad 0\leq t\leq T ,
\end{equation*}
where $\vartheta =\left(\alpha ,\beta \right)\in \Theta
=\left(a_1,a_2\right)\times \left(b_1,b_2\right)$, $b_1>0$ and $\gamma
\geq {1}$.

Recall that the MLE $\hat\vartheta
_T=\left(\hat\alpha _T,\hat\beta _T\right)$ of the parameter $\vartheta $ is
consistent and asymptotically normal. Moreover, the moments of this estimator
converge too (see Theorem 2.8, \cite{Kut04}):
\begin{equation}
\label{10}
\sup_{\theta\in \Theta}  T^{\frac{p}{2}}\Ex_\vartheta \left|\hat\vartheta _T-\vartheta \right|^p\leq C,
\end{equation}
for any $p>0$.

\subsubsection{The Test Based on Empirical Distribution Function}

We study the test
$$
\hat\psi _T\left(X^T\right)=\1_{\left\{\Delta _T\left(X^T\right)>c_\varepsilon \right\}},
$$
where the  test statistic is
\begin{equation*}
\label{Del}
\Delta _T\left(X^T\right)=\hat\beta_{T}^{\frac{2}{\gamma +1}}
  \sigma^{\frac{2\left(\gamma -1\right)}{\gamma
      +1}}T\int_{-\infty }^{\infty }\left[\hat
  F_T\left(x\right)-F\left(\hat \vartheta _T,x\right)\right]^2{\rm
  d}F\left(\hat \vartheta _T,x\right) .
\end{equation*}

Let us introduce the random variable
\begin{align*}
\Delta =\int_{-\infty }^{\infty }\left[ \Phi \left(y\right)+
  \frac{\Pi}{\gamma a}\; f_0\left(y\right)
 +
    \frac{y\;\Psi}{\left(\gamma +1\right) b}\; f_0\left(y\right)  \right]^2
f_0\left(y\right){\rm d}y,
\end{align*}
where
\begin{align*}
\Phi\left(y\right)&=2\int_{-\infty }^{\infty
}\frac{F_0\left(z\right)F_0\left(y\right)-
F_0\left(z\wedge y\right)}{\sqrt{f_0\left(y\right)}}\;{\rm d}W\left(z\right),\\
 \Pi&=\int_{-\infty }^{\infty
}\left|z\right|^{\gamma-1} \sqrt{f_0\left(z\right)}\;{\rm
  d}W\left(z\right),\qquad \quad a=\Ex_0\left|\xi\right|
  ^{2\gamma-2 },\\
 \Psi&=\int_{-\infty }^{\infty }
 \sgn\left(z\right)\left|z\right|^\gamma \sqrt{f_0\left(z\right)}\;{\rm
   d}W\left(z\right),\quad b=\Ex_0\left|\xi\right|
  ^{2\gamma } .
\end{align*}
Here $W\left(\cdot \right)$ is two-sided Wiener process.  The constant
$c_\varepsilon $ is defined by the equation
$$
\Pb\left(\Delta >c_\varepsilon \right)=\varepsilon .
$$
The distribution of the random variable $\Delta $ is not known in a closed form but the value
$c_\varepsilon $ can be easily obtained with the help of the Monte Carlo
simulations. Let us stress that this value is the same for all $\vartheta $
and therefore can be calculated before the experiment.

Our first result is

\begin{theorem}
\label{T1} The test $\hat\psi _T\left(X^T\right)\in {\cal K}_\varepsilon $.
\end{theorem}
{\bf Proof.} We have the relation
\begin{align*}
&\sqrt{T}\left(\hat F_T\left(x\right)-F\left(\hat \vartheta _T,x\right)
\right)\\
&\qquad =\sqrt{T}\left(\hat F_T\left(x\right)-F\left( \vartheta ,x\right)
\right)+\sqrt{T}\left(F\left( \vartheta ,x\right) -F\left(\hat \vartheta
_T,x\right) \right)\\
&\qquad =\eta _T\left(x\right)-\left(\sqrt{T}\left(\hat\vartheta_T- \vartheta
 \right),\frac{\partial F\left( \vartheta ,x\right)}{\partial \vartheta }\right) +r_T.
\end{align*}
Here $\left(\sqrt{T}(\hat\vartheta_T- \vartheta
 ),\frac{\partial F\left( \vartheta ,x\right)}{\partial \vartheta
 }\right) $ is the usual scalar product in $R^2$ and stochastic process $\eta
 _T\left(x\right) $ is given by
$$
\eta _T\left(x\right)=\frac{2}{\sqrt{T}}\int_{0}^{T }\frac{F\left(\vartheta ,X_t\right)F\left(\vartheta ,x\right) -F\left(\vartheta ,X_t\wedge x\right)}{\sigma
  \,f\left(\vartheta ,X_t\right)}\,{\rm d}W_t.
$$
The convergence $r_T\rightarrow 0$ follows from the  representation
\eqref{df}, the estimate
\begin{equation}
\label{x}
\sup_{\vartheta \in\Theta} \Ex_\vartheta \left(\int_{0}^{\xi
}\frac{F\left(\vartheta ,y\right)F\left(\vartheta ,x\right) -F\left(\vartheta
  ,y\wedge x\right)}{
  \,f\left(\vartheta ,y\right)}\,{\rm d}y\right)^2<C,
\end{equation}
which can be obtained by direct calculation (see \cite{Kut04}, Theorem 4.6)
and the estimate
\eqref{10}.
Note that the density
$f\left(\vartheta ,x\right)$ has exponentially decreasing tails and all
necessary estimates can be derived in the straightforward way.

Define the random functions
\begin{align*}
\hat\eta _T\left(x\right)&= \hat\beta_{T}^{\frac{1}{\gamma +1}}
\sigma^{\frac{\gamma -1}{\gamma +1}}\left[\eta
  _T\left(x\right)-\left(\sqrt{T}\left(\hat\vartheta_T- \vartheta
  \right),\frac{\partial F\left( \vartheta ,x\right)}{\partial \vartheta
  }\right)\right],\qquad x\in R,\\ \eta_0 \left(x\right)&=\Phi \left(y\right)+
  \frac{\Pi}{\gamma a}\; f_0\left(y\right)
 +
    \frac{y\;\Psi}{\left(\gamma +1\right) b}\; f_0\left(y\right) ,\qquad x\in R.
\end{align*}
We have to verify that
\begin{equation}
\label{11}
\int_{-\infty }^{\infty }\hat\eta _T\left(x\right)^2f(\hat
\vartheta _T,x)\;{\rm d}x\Longrightarrow \int_{-\infty }^{\infty }\eta
_0\left(x\right)^2f_0\left(x\right)\,{\rm d}x.
\end{equation}

We start with the convergence of the finite dimensional distributions.
The form of the invariant density \eqref{pdf} suggests the following change of variables
$$
Z_t=\frac{\beta ^{\frac{1}{\gamma +1}}}{\sigma ^{\frac{2}{\gamma +1}}}\left(X_t-\alpha \right).
$$
This process satisfies the following stochastic differential
$$
{\rm d}Z_t=-\sgn\left(Z_t\right)\;   \left|Z_t\right|^{\gamma }\;{\rm d}\left(t\beta ^{\frac{2}{\gamma +1}}\sigma
^{\frac{2\left(\gamma -1\right)}{\gamma +1}}\right) +\beta ^{\frac{1}{\gamma +1}}\sigma
^{\frac{\gamma -1}{\gamma +1}}{\rm d}W_t,\quad 0\leq t\leq T.
$$
Therefore, if we denote
$$
Y_s=Z_{s\beta ^{-\frac{2}{ \gamma +1}}\sigma
  ^{-\frac{2\left(\gamma -1\right)}{\gamma +1}}},\qquad 0\leq s=t\beta ^{\frac{2 }{\gamma +1}}\sigma
^{\frac{2\left(\gamma -1\right)}{\gamma +1}}\leq T_*=T\beta ^{\frac{2 }{\gamma +1}}\sigma
^{\frac{2\left(\gamma -1\right)}{\gamma +1}},
$$
the process $Y_s$ satisfies the equation
$$
{\rm d}Y_s=-\sgn\left(Y_s\right)\;\left|Y_s\right|^{\gamma }{\rm d}s+{\rm d}w_s,\qquad Y_0,\quad 0\leq s\leq T_*,
$$
where $w_s=\beta ^{\frac{1}{\gamma +1}}\sigma
^{\frac{\gamma -1}{\gamma +1}}W_t$ (here $t =s\beta ^{-\frac{2 }{\gamma +1}}\sigma
^{-\frac{2\left(\gamma -1\right)}{\gamma +1}} $) is another Wiener process.
Obviously the
 process $Y_s$  is ergodic with the  invariant density
$
f_0\left(x\right).
$

Let us define  $y=\beta ^{\frac{1}{\gamma +1}}\sigma ^{-\frac{2}{\gamma +1}} \left(x-\alpha
\right)$. Then we can write
\begin{align*}
&F\left(\vartheta ,x\right)=F_0\left(y\right),\qquad \quad f\left(\vartheta
,x\right)=\beta ^{\frac{1}{\gamma +1}}\sigma ^{-\frac{2}{\gamma +1}}
f_0\left(y\right).\\
\end{align*}
For the stochastic process $\eta _T\left(x\right) $ this change of variables
gives the representation
\begin{align*}
\eta _T\left(x\right)&=\frac{2}{\sqrt{T}}\int_{0}^{T}\frac{F\left(\vartheta
  ,X_t\right)F\left(\vartheta ,x\right)-F\left(\vartheta ,x\wedge
  X_t\right)}{\sigma f\left(\vartheta ,X_t\right)} \;{\rm
  d}W_t\\ &=\frac{2}{\sqrt{T}}\int_{0}^{T}\frac{F_0\left(Z_t\right)F_0\left(y\right)-F_0\left(y\wedge
  Z_t\right)}{\beta^{\frac{1}{\gamma +1}} \sigma^{\frac{\gamma -1}{\gamma +1}}
  f_0\left(Z_t\right)} \;{\rm
  d}W_t\\ &=\frac{2\beta^{-\frac{1}{\gamma +1}}
  \sigma^{-\frac{\gamma -1}{\gamma
      +1}}}{\sqrt{T_*}}\int_{0}^{T_*}\frac{F_0\left(Y_s\right)F_0\left(y\right)-F_0\left(y\wedge
  Y_s\right)}{
  f_0\left(Y_s\right)} \;{\rm d}w_s\\ &=\beta^{-\frac{1}{\gamma +1}}
  \sigma^{-\frac{\gamma -1}{\gamma
      +1}} \;\Phi _{T_*}\left(y\right) ,
\end{align*}
where the last equality defines the random function $\Phi
_{T_*}\left(y\right) $.

Introduce the following integrals
\begin{align*}
&\pi  _{T_*}=\frac{1}{\sqrt{T_*}}\int_{0}^{T_*}\left|Y_s\right|^{\gamma
  -1}{\rm d}w_s ,\quad \psi
  _{T_*}=\frac{1}{\sqrt{T_*}}\int_{0}^{T_*}\sgn\left(Y_s\right)\,\left|Y_s\right|^{\gamma
  }{\rm d}w_s  ,\\
&a_{T_*}=\frac{1}{{T_*}}\int_{0}^{T_*}\left|Y_s\right|^{2\gamma -2
}{\rm d}s,\quad\quad  b_{T_*}=\frac{1}{{T_*}}\int_{0}^{T_*}\left|Y_s\right|^{2\gamma
}{\rm d}s,\\
&c_{T_*}=\frac{1}{{T_*}}\int_{0}^{T_*}\sgn\left(Y_s\right)\left|Y_s\right|^{2\gamma -1
}{\rm d}s.
\end{align*}
Note that by the law of large numbers we have
$$
a_{T_*}\longrightarrow a,\qquad \quad b_{T_*}\longrightarrow b
$$
The invariant density $f_0\left(y\right)$ is a symmetric function and therefore
\begin{equation}
\label{n}
c_{T_*}\longrightarrow \Ex_{\vartheta _0} \left(\sgn\left(\xi \right)
\left|\xi \right|^{2\gamma -1}\right)=0.
\end{equation}
The random variables $\pi  _{T_*} $ and $\psi  _{T_*}$  are asymptotically normal
by the central limit theorem
$$
\pi  _{T_*}\Longrightarrow \Pi ,\quad \qquad \psi
_{T_*}\Longrightarrow \Psi
$$
and  due to \eqref{n} they are  asymptotically independent.

 The MLE $\hat\vartheta  _T$ admits the following representation
\begin{equation}
\label{mle}
\sqrt{T}\left(\hat\vartheta_T- \vartheta \right)={\rm I}_T\left(\vartheta
\right)^{-1}\frac{1 }{\sigma\sqrt{T}}\int_{0}^{T}\frac{\partial S\left(\vartheta
  ,X_t\right)}{\partial \vartheta
}\;{\rm d}W_t+o\left(1\right),
\end{equation}
 where ${\rm I}_T\left(\vartheta \right) $ is the
  $2\times 2$  matrix
$$
{\rm I}_T\left(\vartheta \right)=\frac{1}{T\sigma ^2}\int_{0}^{T}\frac{\partial S\left(\vartheta
  ,X_t\right)}{\partial \vartheta  }\left(\frac{\partial S\left(\vartheta
  ,X_t\right)}{\partial \vartheta  }\right)^\tau\;{\rm d} t.
$$
Here $\tau $ means transposition. For  the proof of this representation see
 Theorem 2.8 in  \cite{Kut04} and Theorem 8.1 in \cite{IH81}. The convergence
 \eqref{n} allows  us to consider the information matrix as asymptotically
 diagonal.

Note that for the trend coefficient $S\left(\vartheta ,x\right)=-\beta\;
\sgn\left(x-\alpha \right)\left|x-\alpha \right|^\gamma $ we have the equality
$
\frac{\partial S\left(\vartheta ,x\right)}{\partial \alpha }=\beta \gamma
\left|x-\alpha \right|^{\gamma -1} .
$
We have
\begin{align*}
&\frac{1}{T\sigma ^2}\int_{0}^{T}\left(\frac{\partial S\left(\vartheta
  ,X_t\right)}{\partial \alpha   }\right)^2\;{\rm d} t=\gamma ^2\,\sigma
  ^{\frac{2\left(\gamma -3\right)}{\gamma +1}}\,\beta ^{\frac{4}{\gamma
      +1}}\,a_{T_*},\\
&\frac{1}{T\sigma ^2}\int_{0}^{T}\left(\frac{\partial S\left(\vartheta
  ,X_t\right)}{\partial \beta   }\right)^2\;{\rm d} t=\sigma
  ^{\frac{2\left(\gamma -1\right)}{\gamma +1}}\,\beta ^{-\frac{2\gamma
    }{\gamma +1}}\,b_{T_*}
\end{align*}
and
\begin{align*}
&\frac{1}{\sqrt{T}\sigma }\int_{0}^{T}\frac{\partial S\left(\vartheta
  ,X_t\right)}{\partial \alpha   }\;{\rm d} W_t=\gamma \,\sigma
  ^{\frac{\gamma -3}{\gamma +1}}\,\beta ^{\frac{2}{\gamma
      +1}}\,\pi _{T_*},\\
&\frac{1}{\sqrt{T}\sigma }\int_{0}^{T}\frac{\partial S\left(\vartheta
  ,X_t\right)}{\partial \beta   }\;{\rm d} W_t=-\sigma
  ^{\frac{\gamma -1}{\gamma +1}}\,\beta ^{-\frac{\gamma }{\gamma
      +1}}\,\psi _{T_*}
\end{align*}
Further
$$
\frac{\partial
  F\left(\vartheta ,x\right)}{\partial \alpha  }=-\frac{\beta
  ^{\frac{1}{\gamma +1}}}{\sigma ^{\frac{2}{\gamma +1}}    }\;f_0\left(y\right) ,\qquad \frac{\partial F\left(\vartheta ,x\right)}{\partial \beta }=\frac{y}{\beta
  \left(\gamma +1\right)}\;f_0\left(y\right) .
$$
Therefore we can write
\begin{align*}
&\left(\sqrt{T}\left(\hat\vartheta  _T-\vartheta\right),\;\frac{\partial F\left(\vartheta
  ,x\right)}{\partial \vartheta}\right)\\
&\quad = \sqrt{T}\left(\hat\alpha   _T-\alpha \right)\;\frac{\partial F\left(\vartheta
  ,x\right)}{\partial \alpha }+\sqrt{T}\left(\hat\beta   _T-\beta \right)\;\frac{\partial F\left(\vartheta
  ,x\right)}{\partial \beta }\\
&\quad =-\beta ^{-\frac{1}{\gamma +1}}  \sigma ^{\frac{1-\gamma }{\gamma +1}}
  \left[\frac{\pi _{T_*}}{\gamma\;  a_{T_*}}+ \frac{y\,\psi _{T_*}}{\left(\gamma+1\right)\;  b_{T_*}} \right]f_0\left(y\right) +o\left(1\right).
\end{align*}

This, in turn, allows  us to write
\begin{align*}
&\sqrt{T}\left(\hat F_T\left(x\right)-F\left(\hat \vartheta _T,x\right)
\right)\\
&\qquad  =\beta^{-\frac{1}{\gamma +1}}
  \sigma^{\frac{1-\gamma }{\gamma
      +1}} \;\left[ \Phi _{T_*}\left(y\right)+\frac{\pi _{T_*}}{\gamma\;  a_{T_*}}f_0\left(y\right)+ \frac{y\,\psi _{T_*}}{\left(\gamma+1\right)\;  b_{T_*}}f_0\left(y\right)
    \right] +o\left(1\right).
\end{align*}
Finally we obtain
\begin{align*}
&\Delta _T\left(X^T\right)\\
&\quad =\int_{-\infty }^{\infty } \left[\Phi _{T_*}\left(y\right)+\frac{\pi
      _{T_*}}{\gamma\;  a_{T_*}}f_0\left(y\right)+ \frac{y\,\psi
      _{T_*}}{\left(\gamma+1\right)\;  b_{T_*}}f_0\left(y\right)
    \right]^2 f_0\left(y\right){\rm d}y+o\left(1\right).
\end{align*}
Now we can replace  $a_{T_*} $ and $b_{T_*}$ with their limits $a$ and $b$
and  denote
$$
\tilde \eta _{T_*}\left(y\right)=\Phi _{T_*}\left(y\right)+\frac{\pi
      _{T_*}}{\gamma\;  a}f_0\left(y\right)+ \frac{y\,\psi
      _{T_*}}{\left(\gamma+1\right)\, b}f_0\left(y\right) .
$$
We shall verify the convergence
\begin{equation}
\label{c}
\int_{-\infty }^{\infty } \tilde \eta
_{T_*}\left(y\right)^2f_0\left(y\right)\,{\rm d}y\Longrightarrow
\int_{-\infty }^{\infty }  \eta
_{0}\left(y\right)^2f_0\left(y\right)\,{\rm d}y.
\end{equation}

To prove it,  we shall check the following three conditions.
\begin{enumerate}

\item  {\it The finite dimensional distributions of $\tilde\eta _{T_*}\left(\cdot
  \right)$ converge, i.e., for any $k\geq 1$ and any $y_1,\ldots,y_k$ we have
  the convergence
$$
\left(\tilde \eta _{T_*}\left(y_1\right),\ldots, \tilde\eta
_{T_*}\left(y_k\right)\right)\Longrightarrow \left( \eta
_0\left(y_1\right),\ldots, \eta _0\left(y_k\right)\right).
$$
\item   There exists a constant $C_1>0$ such
that
\begin{align}
\label{12}
&\Ex_{\vartheta _0}\left|\tilde\eta _{T_*}\left(y_2\right)-\tilde\eta
_{T_*} \left(y_1\right)\right|^2\leq C_1\,\left|y_2-y_1\right| .
\end{align}

\item There exist constants $ C_2>0$ and $\kappa >0$ such }
\begin{align}
\label{13}
&\Ex_{\vartheta _0}\left|\tilde\eta _{T_*}\left(y\right)\right|^2\leq {C}_2\,e^{-\kappa \left|y\right|^{\gamma +1}}.
\end{align}
\end{enumerate}

 If these conditions hold,  \eqref{c}  follows  from the
 results of \cite{IH81}. Indeed, by
 Theorem A.22 in \cite{IH81}   integrals converge on any
 finite interval $\left[-L,L\right]$ and outside of this interval we can
 estimate the  tail integrals   as   in
 the proof of the Theorem  1.5.6  \cite{IH81}.

By the central limit theorem for stochastic integral we obtain the desired
joint asymptotic normality
\begin{align*}
&\left(\Phi _{T_*}\left(y_1\right),\ldots , \Phi
  _{T_*}\left(y_k\right), \pi _{T_*},\psi _{T_*} \right)\Longrightarrow
  \left(\Phi\left(y_1\right),\ldots , \Phi \left(y_k\right),\Pi ,\Psi \right),
\end{align*}
which proves  convergence of the finite dimensional distributions of
$\tilde\eta _T\left(\cdot \right)$.
Further, for all $y_1,y_2, \left|y_2-y_1\right|\leq 1$  we have the
estimate   ($y_1<y_2$)
\begin{align*}
&\Ex_{\vartheta _0} \left|\tilde\eta _T\left(y_2 \right)-\tilde\eta _T\left(y_1
\right) \right|^2 \leq 3\Ex_{\vartheta _0} \left|\Phi _{T_*}\left(y_2\right)-\Phi
_{T_*}\left(y_1\right)\right| ^2\\
&\qquad \quad
+\frac{3}{\left(\gamma +1\right)^2b^2}\left|y_2f_0\left(y_2\right)-y_1f_0\left(y_1\right)   \right| ^2\Ex_{\vartheta _0}
\psi _{T_*}^2\\
&\qquad \quad
+\frac{3}{\gamma ^2a^2}\left|f_0\left(y_2\right)-f_0\left(y_1\right)\right| ^2\Ex_{\vartheta _0}
\pi _{T_*}^2  \leq C\left|y_2-y_1\right|^2
\end{align*}
because
\begin{align*}
&\Ex_{\vartheta _0}\left|\Phi _{T_*}\left(y_2\right)-\Phi
  _{T_*}\left(y_1\right)\right| ^2=4\int_{-\infty
  }^{y_1}\frac{F_0\left(z\right)^2\left[F_0\left(y_2\right)-F_0\left(y_1\right)\right]^2}{f_0\left(z\right)}{\rm
    d}z\\ &\qquad \quad+ 4\int_{y_2 }^{\infty
  }\frac{F_0\left(z\right)^2\left[F_0\left(y_2\right)-F_0\left(y_1\right)\right]^2}{f_0\left(z\right)}{\rm
    d}z\\ &\qquad \quad+ 4\int_{y_1 }^{y_2
  }\frac{\left(F_0\left(z\right)\left[F_0\left(y_2\right)-F_0\left(y_1\right)\right]-F_0\left(z\right)+F_0\left(y_1\right)\right)^2}{f_0\left(z\right)}{\rm
    d}z\\ &\qquad\leq C\left|y_2-y_1\right|^2.
\end{align*}
\bigskip
Therefore, we obtain  \eqref{12}. To prove  \eqref{13}, write
\begin{align*}
&\Ex_{\vartheta _0} \tilde\eta _T\left(y\right)^2\leq 3   \Ex_{\vartheta _0}\Phi
_{T_*}\left(y\right) ^2 +\frac{3}{\left(\gamma +1\right)^2b^2} y ^2f_0\left(y\right) ^2\Ex_{\vartheta _0}
\Psi _{T_*}^2+\frac{3}{\gamma ^2a^2} f_0\left(y\right) ^2\Ex_{\vartheta _0}
\Pi _{T_*}^2\\
&\quad \leq 12\int_{-\infty }^{\infty
}\frac{\left[F_0\left(z\right)F_0\left(y\right)-F_0\left(z\wedge
    y\right)\right]^2}{f_0\left(z\right) }{\rm d}z  +C\left
(1+y^2\right)f_0\left(y\right) ^2 .
\end{align*}
Further, using the same arguments as in \cite{Kut04}, Example 4.1.3, we obtain the
estimate
$$
\left(F_0\left(y\right)-1\right)^2\int_{-\infty }^{y
}\frac{F_0\left(z\right)^2}{f_0\left(z\right)
}{\rm d}z  +F_0\left(y\right)^2\int_{y }^{\infty
}\frac{\left[F_0\left(z\right)-1\right]^2}{f_0\left(z\right) }{\rm d}z\leq
C\,e^{-\kappa  \left|y\right|^{\gamma +1}}
$$
with some constants $C>0,\kappa >0$. For example, for the large values of $z$ we can
write
\begin{align*}
&\frac{\left[F_0\left(z\right)-1\right]^2}{f_0\left(z\right) }=\left(\int_{z}^{\infty
} \exp\left\{-cu^{\gamma +1}+\frac{c}{2}z^{\gamma +1}\right\}{\rm d}u\right)^2
  \\
&\qquad \leq \left(\int_{z}^{\infty
} \exp\left\{-\frac{c}{2}u^{\gamma +1}\right\}{\rm d}u\right)^2\leq C\,e^{-cz^{\gamma +1}}
\end{align*}
and so on.

Therefore the conditions of the weak convergence of integrals are verified
and the test $\hat\psi_T\in{\cal K}_\varepsilon $.

The consistency of this test is implied by the elementary inequalities as
follows. Suppose that the trend coefficient $S\left(x\right)$ of the observed process \eqref{7}
does not belong to the given  parametric family
$S\left(\vartheta ,x\right),\vartheta \in  \left[a_1,a_2\right]\times
\left[b_1,b_2\right] $. It is known (Proposition 2.36, \cite{Kut04}) that the
MLE $\hat\vartheta _T$ converges
to the value $\hat\vartheta $ which minimizes the Kullback-Leibner distance
between the parametric family and the true distribution:
$$
\hat\vartheta =\arg\inf_{\vartheta }\Ex_S\left(\frac{S\left(\xi_*
  \right)-S\left(\vartheta ,\xi_*  \right)}{\sigma \left(\xi*
  \right)}\right)^2 .
$$
Here the random variable  $\xi _*$ has the invariant density function
$f_S\left(x\right)$.  It can be shown that
$$
\int_{R}^{ }\left[\hat
  F_T\left(x\right)-F(\hat \vartheta _T,x)\right]^2{\rm
  d}F(\hat \vartheta _T,x) \rightarrow \int_{R }^{ }\left[
  F\left(x\right)-F(\hat \vartheta ,x)\right]^2{\rm
  d}F(\hat \vartheta,x)>0.
$$
Therefore for  the statistic $\Delta _T$ we have
$$
\Delta _T\longrightarrow \infty \qquad {\rm and}\qquad \Pb_S\left\{\Delta
_T>c_\varepsilon \right\}\rightarrow 1.
$$
Hence the test is consistent.

\subsubsection{The Test Based on Empirical Density}

Now we study the test
$$
\tilde\psi _T\left(X^T\right)=\1_{\left\{\delta _T\left(X^T\right)>c_\varepsilon \right\}},
$$
where the  test statistic is
\begin{equation*}
\label{del}
\delta _T\left(X^T\right)=
  \sigma^{2}T\int_{-\infty }^{\infty }\left[\hat
  f_T\left(x\right)-f\left(\hat \vartheta _T,x\right)\right]^2{\rm
  d}F\left(\hat \vartheta _T,x\right) .
\end{equation*}
Let us denote $\zeta _T\left(\vartheta ,x\right)=\sqrt{T}\left(\hat
  f_T\left(x\right)-f\left( \vartheta
  ,x\right)\right)$ and write
\begin{align*}
\sqrt{T}\left(\hat  f_T\left(x\right)-f\left(\hat \vartheta _T,x\right)\right)&=\zeta
_T\left(\vartheta ,x\right) -\sqrt{T}\left(\hat\alpha _T-\alpha
\right)\frac{\partial f\left(\vartheta ,x\right)}{\partial \alpha  }\\
&\quad -\sqrt{T}\left(\hat\beta  _T-\beta  \right)\frac{\partial
  f\left(\vartheta ,x\right)}{\partial \beta   } +o\left(1\right).
\end{align*}

  Using the same arguments as above we obtain the representations
\begin{align*}
\zeta _T\left(\vartheta ,x\right)&=\frac{2f_0\left(y\right)}{\sigma
    \sqrt{T_*}}\int_{0}^{T_*}\frac{F_0\left(Y_s\right)-\1_{\left\{Y_s>y\right\}}}{f_0\left(Y_s\right)}{\rm
    d}w_s+o\left(1\right)\\
&=\sigma ^{-1}\tilde\Phi_{T_*}\left(y\right)f_0\left(y\right)+o\left(1\right) ,\\
\frac{\partial f\left(\vartheta ,x\right)}{\partial \alpha
}&=-2\beta^{\frac{2}{\gamma +1}}\;\sigma ^{-\frac{2}{\gamma +1}}
\,\sgn\left(y\right)\,\left|y\right|^{\gamma }f_0\left(y\right),\\
\frac{\partial f\left(\vartheta ,x\right)}{\partial \beta }&=-\frac{\beta^{-\frac{\gamma }{\gamma +1}}\sigma ^{-\frac{2}{\gamma +1}} }{\gamma +1}\left[1-2\left|y\right|^{\gamma +1}\right]f_0\left(y\right) .
\end{align*}
These equalities together with the representations   of the estimators
$\hat\alpha _T$ and $\hat \beta _T$ allow us to  write
\begin{align*}
&\sigma\sqrt{T}\left(\hat  f_T\left(x\right)-f(\hat \vartheta
_T,x)\right)\\
&\qquad = \left[\tilde\Phi
_{T_*}\left(y\right)+2 \sgn\left(y\right)
\left|y\right|^{\gamma }\frac{\pi_{T_*}}{\gamma a}
+\left[1-2\left|y\right|^{\gamma
    +1}\right]\frac{\psi _{T_*}}{\left(\gamma +1\right)b}\right]f_0\left(y\right) +o\left(1\right).
\end{align*}
Using the same arguments as in the section 3.1.1 the following convergence
\begin{equation*}
\label{ltdel}
\delta _T\Longrightarrow \delta =\int_{-\infty }^{\infty }\zeta
_0\left(y\right)^2\,f_0\left(y\right)\,{\rm d}y  ,
\end{equation*}
can be proved. Here
\begin{align*}
\zeta _0\left(y\right)&=\left[\tilde\Phi
\left(y\right)+2 \sgn\left(y\right)
\left|y\right|^{\gamma }\frac{\Pi}{\gamma a}+\left[1-2\left|y\right|^{\gamma
    +1}\right]\frac{\Psi}{\left(\gamma +1\right)b}
  \right]f_0\left(y\right),\\
\tilde\Phi\left(y\right)&=2\int_{-\infty }^{\infty
}\frac{F_0\left(z\right)-\1_{\left\{z>y \right\}}}{ \sqrt{f_0\left(z\right)}}
\,{\rm d}W\left(z\right).
\end{align*}
Hence the test $\tilde \psi _T=\1_{\left\{\delta _T>c_\varepsilon
  \right\}}$ is APF. The threshold $c_\varepsilon $ is defined by  the equation
$$
\Pb\left\{\delta >c_\varepsilon \right\}=\varepsilon
$$
and  therefore  it belongs to ${\cal
  K}_\varepsilon $.

\subsection{Case $0\leq \gamma <\frac{1}{2}$.}
If we observe (under hypothesis ${\cal H}_0$) the same equation
$$
{\rm d}X_t=-\beta \,\sgn\left(X_t-\alpha \right)\left|X_t-\alpha\right|^\gamma
\,{\rm d}t+\sigma \,{\rm
  d}W_t,\quad X_0,\quad 0\leq t\leq T .
$$
but with $\gamma \in [0,\frac{1}{2})$, then the main difference with $\gamma
  \geq 1$ is due to the rate of convergence of the MLE $\hat\alpha _T$. As the
  rate is faster than $\sqrt{T}$ the contribution of this estimator to the
  limit distribution of test statistic is negligeable. This property of the
  test statisics
  was mentioned by Darling \cite{Dar}.

Recall that  in the
  case $\gamma =0$ (simple switching)
$$
{\rm d}X_t=-\beta \,\sgn\left(X_t-\alpha \right)\,{\rm d}t+\sigma \,{\rm
  d}W_t,\quad X_0,\quad 0\leq t\leq T .
$$
we have the convergence
$$
T\left(\hat\alpha  _T-\alpha  \right)\Longrightarrow \hat u,\qquad
Z\left(\hat u\right)=\sup_uZ\left(u\right)
$$
where
$$
Z\left(u\right)=\exp\left\{W\left(u\gamma _\vartheta
\right)-\frac{\left|u\gamma _\vartheta\right|}{2}\right\}.
$$
Here $W\left(\cdot \right)$ is double sided Wiener process and  $\gamma
_\vartheta >0$ is some constant (see details in \cite{Kut04}, Section 3.4).
Moreover, we have the convergence of moments too: for any $p>0$
$$
\Ex_\vartheta \left|T\left(\hat\alpha  _T-\alpha
\right)\right|^p\longrightarrow  \Ex_\vartheta\left|\hat u\right|^p.
$$
Therefore, if we repeat the proofs  above,  we shall see that
$$
\sqrt{T}\left(\hat\alpha  _T-\alpha
\right)\frac{\partial F\left(\vartheta ,x\right)}{\partial \alpha
}=O\left(\frac{1}{\sqrt{T}}\right)
$$
and
$$
 \sqrt{T}\left(\hat\alpha  _T-\alpha
\right)\frac{\partial f\left(\vartheta ,x\right)}{\partial \alpha
}=O\left(\frac{1}{\sqrt{T}}\right).
$$
Of course, we have to be careful with the second term because the invariant
density  is
\begin{equation*}
f\left(\vartheta ,x\right)=\frac{\beta}{\sigma ^2}\, \exp\left\{{-\frac{2\beta}{\sigma
    ^2} \left|x-\alpha \right|} \right\}
\end{equation*}
and the derivative is not continuous. However, the function $f\left(\vartheta
,x\right)$ is absolutely continuous and this is sufficient for the proof.

The  Cram\'er-von Mises type statistics are
 \begin{align*}
\Delta _T\left(X^T\right)&= \frac{\hat\beta _T^2}{\sigma ^2}\;T\int_{-\infty }^{\infty
 }\left[\hat F_T\left(x\right)-F\left(\hat\vartheta _T,x\right)\right]^2{\rm
   d} F\left(\hat\vartheta _T,x\right),\\
\delta _T\left(X^T\right)&= {\sigma ^2}\;{T}\int_{-\infty }^{\infty
 }\left[\hat f_T\left(x\right)-f\left(\hat\vartheta _T,x\right)\right]^2{\rm
   d} F\left(\hat\vartheta _T,x\right),
\end{align*}
and their limits are
\begin{align*}
\Delta _T\left(X^T\right)&\Longrightarrow \int_{-\infty }^{\infty }\left[ \Phi \left(y\right)+
    y\;f_0\left(y\right)\;\frac{\Psi}{\left(\gamma +1\right)b}
    \right]^2 f_0\left(y\right){\rm d}y,\\
 \delta _T\left(X^T\right)&\Longrightarrow \int_{-\infty }^{\infty }\left[\tilde\Phi
\left(y\right)+\left[1-2\left|y\right|\right]\frac{\Psi}{\left(\gamma +1\right)b}
  \right]^2f_0\left(y\right)^3{\rm d}y.
\end{align*}
Hence the corresponding tests belong to ${\cal K}_\varepsilon $.

If $\gamma \in\left(0,\frac{1}{2}\right)$, then we have a different limit
$$
T^{\frac{1}{2\gamma
    +1}}\left(\hat\alpha _T-\alpha \right)\Longrightarrow  \hat u,\qquad
Z\left(\hat u\right)=\sup_uZ\left(u\right)
$$
with
$
Z\left(u\right)=\exp\left\{W^H\left(u\Gamma _\vartheta
\right)-\frac{\left|u\Gamma _\vartheta\right|^{2H}}{2}\right\}.
$
Here $W^H\left(\cdot \right)$ is two sided fractional Brownian motion,
$H=\kappa +\frac{1}{2} $ (Hurst parameter)  and  $\Gamma
_\vartheta$ is some constant. We have the convergence of moments too: for any $p>0$
$$
T^{\frac{p}{2\gamma
    +1}}\Ex_\vartheta \left|\left(\hat\vartheta _T-\vartheta
\right)\right|^p\longrightarrow  \Ex_\vartheta\left|\hat u\right|^p.
$$
For the proofs see \cite{DK03} or \cite{Kut04}, Section 3.2.

Hence once again we can use the tests $\hat\psi _T$ and $\tilde\psi _T$ and
 the limits of the test statistics are obtained by setting $\Pi \equiv 0$ in
 $\Delta $ and $\delta $ respectively.

\subsection{Discussion}

The case   $\gamma \in [\frac{1}{2},1)$ was not included in this study because the
  appropriate properties of the  MLE $\hat \alpha _T$ are available  only in the cases   $\gamma \in
  [0,\frac{1}{2})$ \cite{DK03} and   $\gamma \geq 1$ \cite{Kut04}.  In
    the case $\gamma \in (\frac{1}{2},1)$ the derivative of the trend with
      respect to parameter $\alpha $ is no
      more locally bounded, but the singularity at the point $x=\alpha  $ is
      integrable and the proof presented in \cite{Kut04}, Theorem 2.8   can be
      carried out. Note that the   Fisher information is bounded. Therefore in
      this case we obtain the same result as for $\gamma \geq {1}$.
This is not the case if $\gamma =\frac{1}{2} $ and for this model we need a
special study. Note that the rate of convergence of the MLE is better
than $\sqrt{T}$ and the limit distribution of the test statistic have to be
the same as for $\gamma \in
  (0,\frac{1}{2})$.

Note that our proofs (strengthened up to the weak convergence
in ${\cal C}_0\left(R\right)$)  imply the following  limits for the
Kolmogorov-Smirnov type statistics  
\begin{align*}
\tilde \Delta _T\left(X^T\right)&=\hat\beta_{T}^{\frac{1}{\gamma +1}}
  \sigma^{\frac{\left(\gamma -1\right)}{\gamma
      +1}}\sup_{x}\sqrt{T}\left|\hat F_T\left(x\right)-F\left(\hat\vartheta
_T,x\right)\right|\Longrightarrow \sup_{y}\left|\eta_0 \left(y\right)\right|,\\
 \tilde \delta _T\left(X^T\right)&={\sigma
}\sup_{x}\sqrt{T}\left|\hat f_T\left(x\right)-f\left(\hat\vartheta
 _T,x\right)\right|\Longrightarrow \sup_{y}\left|\zeta_0 \left(y\right)\right|
\end{align*}
where the limit distributions do  not depend  on $\vartheta $.
 Hence the goodness of fit tests based on these statistics are APF.
It will be interesting to consider  other models with APF tests.

\section{Asymptotically Distribution Free Tests}

Now we consider the second statement of the goodness of fit hypotheses testing
problem. The basic
hypothesis ${\cal H}_0$ is simple: the observed diffusion process satisfies
the stochastic differential equation
$$
{\rm d}X_t=S_0\left(X_t\right)\:{\rm d}t+\sigma \left(X_t\right)\:{\rm
  d}W_t,\quad X_0,\quad 0\leq t\leq T ,
$$
where $S_0\left(\cdot \right)$ and $\sigma \left(\cdot \right)$ are known
functions. We assume that the conditions ${\cal ES}$  and ${\cal A}_0$
hold. Therefore, the stochastic process are ergodic   with
the invariant density 
$$
f_{S_0}\left(x\right)=\frac{1}{G\left(S_0\right)\sigma
  \left(x\right)^2}\;\exp\left\{2\int_{0}^{x} \frac{S_0\left(y\right)}{\sigma
  \left(y\right)^2}\; {\rm d}y\right\}.
$$
Of course, we can use the Cram\'er-von Mises type statistics based on the empirical density
$$
\delta _T\left(X^T\right)=T\int_{-\infty }^{\infty }\left[\hat
  f_T\left(x\right)-f_{S_0}\left(x\right)\right]^2 {\rm d}F_{S_0}\left(x\right),
$$
but its limit under hypothesis is
$$
\delta \left(S_0\right)=\int_{-\infty }^{\infty }\zeta \left(S_0,x\right)^2 {\rm
  d}F_{S_0}\left(x\right),
$$
where
$$
\zeta
\left(S_0,x\right)=2f_{S_0}\left(x\right) \int_{-\infty }^{\infty
}\frac{F_{S_0}\left(y\right)- \1_{\left\{y>x\right\}}}{\sigma \left(y\right)
  \sqrt{f_{S_0}\left(y\right)}}\; {\rm d}W\left(y\right)
$$
is the limit in distribution of the normalized difference $$\zeta_T
\left(S_0,x\right)= \sqrt{T}\left(\hat
  f_T\left(x\right)-f_{S_0}\left(x\right)\right)$$ (see the representation
  \eqref{lte}).

Therefore, for the test $\hat \phi_T{\left(X^T\right)}=\1_{\left\{ \delta
  _T\left(X^T\right)>c_\varepsilon \right\} } $ we have to solve the equation
$$
\Pb_{S_0}\left(\delta \left(S_0\right)>c_\varepsilon \right)=\varepsilon
$$
and we see that the threshold $c_\varepsilon =c_\varepsilon
\left(S_0\right)$.

Recall that for the i.i.d. observations the limit of the
corresponding statistics based on empirical distribution function is
$$
\Delta _n=n\int_{-\infty }^{\infty }\left[\hat
  F_n\left(x\right)-F_{0}\left(x\right)\right]^2 {\rm
  d}F_{0}\left(x\right)\Longrightarrow  \Delta =\int_{0}^{1}B\left(t\right)^2{\rm d}t ,
$$
where $B\left(t\right), 0\leq t\leq 1$ is the Brownian bridge. Hence the test
$\hat\phi_n=\1_{\left\{\Delta _n>c_\varepsilon \right\}}  $ with $c_\varepsilon $ from
the equation $\Pb\left(\Delta >c_\varepsilon \right)=\varepsilon $ belongs to
the class ${\cal K}_\varepsilon $. Tests based on statistics with limit
distributions independent of the model under hypothesis are called {\it
  asymptotically distribution free} (ADF). There are several works devoted to
the construction of ADF tests for ergodic diffusion processes observed in
continuous time. We can mention here \cite{Fou92},\cite{FouK93},
\cite{Kut04},\cite{DK07}, \cite{NN},\cite{Kut10}, but the connection  of these tests
with the classical Cram\'er-von Mises and Kolmogorov-Smirnov is not evident. One exception is the work \cite{Kut10}, where 
a linear transformation of the normalized deviations $\zeta _T\left(\cdot \right)$ and $\eta  _T\left(\cdot
\right)$    allowed construction of ADF tests. Unfortunately, the
proof in \cite{Kut10} is not satisfactory, since it used a property of the
time change in the Wiener integral, which is not always true. That is why we
decided to suggest     another linear transformation which leads to ADF
test.

Therefore, our goal is to find a linear transformation $L\left(\zeta _T\right)$ of the
random function $\zeta _T\left(x\right)$ such that
$$
\delta _T\left(X^T\right)=\int_{-\infty }^{\infty }\left[L\left(\zeta _T\right)\left(x\right)\right]^2{\rm
  d}F_{S_0}\left(x\right)\Longrightarrow  \int_{0}^{1}w\left(t\right)^2\,{\rm d}t ,
$$
where $w\left(t\right),0\leq t\leq 1 $ is a Wiener process. Then obviously the test
$\hat\psi _T\left(X^T\right)=\1_{\left\{ \delta
  _T\left(X^T\right)>c_\varepsilon \right\}}$ with $c_\varepsilon $ from the
equation
$$
\Pb\left\{\int_{0}^{1}w\left(t\right)^2\,{\rm d}t>c_\varepsilon  \right\}=\varepsilon
$$
will be  ADF.

Let us rewrite the stochastic integral as  follows
\begin{align*}
&\int_{-\infty }^{\infty
}\frac{F_{S_0}\left(y\right)- \1_{\left\{y>x\right\}}}{\sigma \left(y\right)
  \sqrt{f_{S_0}\left(y\right)}}\; {\rm d}W\left(y\right)\\
&\qquad =\int_{-\infty }^{\infty
}\frac{F_{S_0}\left(y\right)- \1_{\left\{F_{S_0}\left(y\right)>F_{S_0}\left(x\right)\right\}}}{\sigma \left(y\right)
  {f_{S_0}\left(y\right)}}\;\sqrt{f_{S_0}\left(y\right)}\; {\rm d}W\left(y\right)\\
&\qquad =\int_{-\infty }^{\infty
}\frac{F_{S_0}\left(y\right)- \1_{\left\{F_{S_0}\left(y\right)>F_{S_0}\left(x\right)\right\}}}{\sigma \left(y\right)
  {f_{S_0}\left(y\right)}}\; {\rm d}w\left(F_{S_0}\left(y\right)\right)\\
&\qquad =\int_{0 }^{1
}\frac{s- \1_{\left\{s>t\right\}}}{a \left(s\right)
  {b\left(s\right)}}\; {\rm d}w\left(s\right)\\
&\qquad =\int_{0 }^{t
}\frac{s}{a \left(s\right)
  {b\left(s\right)}}\; {\rm d}w\left(s\right)+\int_{t}^{1
}\frac{s- 1}{a \left(s\right)
  {b\left(s\right)}}\; {\rm d}w\left(s\right)\equiv u\left(t\right),
\end{align*}
where $w\left(t\right),0\leq t\leq 1$ is a  Wiener process
$$
w\left(t \right)=\int_{-\infty
}^{F_{S_0}^{-1}\left(t\right)}\sqrt{f_{S_0}\left(y\right)}\,{\rm d}W\left(y\right)
$$
and we denoted
\begin{align*}
F_{S_0}\left(y\right)=s,\qquad F_{S_0}\left(x\right)=t,\quad
a\left(s\right)=\sigma \left(F_{S_0}^{-1}\left(s\right) \right),\quad
b\left(s\right)= f_{S_0} \left(F_{S_0}^{-1}\left(s\right) \right).
\end{align*}
Here $F_{S_0}^{-1}\left(s\right) $ is the function inverse to
$F_{S_0}\left(y\right)$, i.e., the solution $y$ of the equation
$F_{S_0}\left(y\right)=s $. Note that under our assumptions the function
$F_{S_0}\left(y\right)$ is strictly increasing.

We can write
$$
{\rm d}u\left(t\right)=\frac{t}{a\left(t\right)b\left(t\right)}\,{\rm
  d}w\left(t\right)-\frac{t-1}{a\left(t\right)b\left(t\right)}\,{\rm
  d}w\left(t\right)=\frac{1}{a\left(t\right)b\left(t\right)}\,{\rm
  d}w\left(t\right) .
$$
Hence the integral (understood in the mean square sense)
$$
\int_{0}^{t} a\left(s\right)b\left(s\right)\;{\rm
  d}u\left(s\right)=w\left(t\right)
$$
provides us the desired transformation. Indeed, we have
$$
u\left(t\right)=\frac{\zeta
  \left(S_0,F_{S_0}^{-1}\left(t\right)
  \right)}{2f_{S_0}\left(F_{S_0}^{-1}\left(t\right)\right)}
$$
and
$$
\int_{0}^{t} a\left(s\right)b\left(s\right)\;{\rm
  d}u\left(s\right)=\int_{-\infty }^{x } \sigma
\left(y\right) f_{S_0}\left(y\right)\;{\rm
  d} \left[\frac{\zeta
  \left(S_0,y
  \right)}{2f_{S_0}\left(y\right)}\right] =w\left(F_{S_0}\left(x\right)  \right)=w\left(t  \right).
$$
Therefore, we can write
\begin{align*}
\delta =\int_{-\infty }^{\infty }\left( \int_{-\infty }^{x } \sigma
\left(y\right) f_{S_0}\left(y\right)\;{\rm
  d} \left[\frac{\zeta
  \left(S_0,y
  \right)}{2f_{S_0}\left(y\right)}\right] \right)^2 {\rm
  d}F_{S_0}\left(x\right)=\int_{0}^{1} w\left(t\right)^2\,{\rm d}t.
\end{align*}

This equality suggests  the statistic
\begin{align*}
\hat \delta_T\left(X^T\right) =\int_{-\infty }^{\infty }\left( \int_{-\infty
}^{x } \sigma \left(y\right) f_{S_0}\left(y\right)\;{\rm
  d} \left[\frac{\zeta_T \left(S_0,y \right)}{2f_{S_0}\left(y\right)}\right]
\right)^2 {\rm d}F_{S_0}\left(x\right),
\end{align*}
where we have to define the integral with respect to the normalized empirical
density
\begin{align*}
\int_{a}^{b} h\left(x\right)\,{\rm d} \left[\frac{\zeta
    _T\left(x\right)}{2f_{S_0}\left(x\right)}\right] .
\end{align*}
If we verify the convergence $\hat \delta_T\left(X^T\right) \Longrightarrow
\delta  $, then  the corresponding test
$$
\hat\psi _T\left(X^T\right)=\1_{\left\{ \hat
  \delta_T\left(X^T\right)>c_\varepsilon  \right\}},\qquad \Pb\left\{ \delta
>c_\varepsilon \right\}=\varepsilon
$$
will be ADF.

Using the representation \eqref{lte} for any piece-wise continuous function
$h\left(x\right)$ with bounded support and partition
$a=x_0<x_1<\ldots<x_m=b$ we can write
\begin{align*}
&\sum_{x_i}^{}h\left(\tilde x_i\right)\left[\frac{\zeta
    _T\left(x_{i+1}\right)}{2f_{S_0}\left(x_{i+1}\right)}-\frac{\zeta
    _T\left(x_{i}\right)}{2f_{S_0}\left(x_{i}\right)}  \right]\\
&\qquad  =\frac{1}{\sqrt{T}}\int_{0}^{T}\frac{\sum_{x_i}^{}h\left(\tilde x_i\right)\1_{\left\{x_i<X_t\leq x_{i+1}\right\}}}{\sigma \left(X_t\right)f_{S_0}\left(X_t\right)}\;{\rm d}W_t\\
&\qquad \qquad \quad
  -\frac{1}{\sqrt{T}}\int_{X_0}^{X_T}\frac{\sum_{x_i}^{}h\left(\tilde
    x_i\right)\1_{\left\{x_i<y\leq x_{i+1}\right\}}}{\sigma
    \left(y\right)f_{S_0}\left(y\right)}\;{\rm d}y.
\end{align*}
Therefore as $\max \left|x_{i+1}-x_i\right|\rightarrow 0$ we obtain the limit
\begin{align*}
&\lim_{ m\rightarrow \infty }\sum_{x_i}^{}h\left(\tilde x_i\right)\left[\frac{\zeta
    _T\left(x_{i+1}\right)}{2f_{S_0}\left(x_{i+1}\right)}-\frac{\zeta
    _T\left(x_{i}\right)}{2f_{S_0}\left(x_{i}\right)}  \right]\\
&\qquad  =\frac{1}{\sqrt{T}}\int_{0}^{T}\frac{h\left(X_t\right)\1_{\left\{a<X_t\leq b\right\}}}{\sigma \left(X_t\right)f_{S_0}\left(X_t\right)}\;{\rm d}W_t
  -\frac{1}{\sqrt{T}}\int_{X_0}^{X_T}\frac{h\left(y\right)\1_{\left\{a<y\leq b\right\}}}{\sigma
    \left(y\right)^2f_{S_0}\left(y\right)}\;{\rm d}y.
\end{align*}
If $a=-\infty $ (as in our case) then this limit exists for a class of functions
vanishing at $-\infty $ because $\sigma
\left(y\right)f_{S_0}\left(y\right)\rightarrow 0$ as
$\left|y\right|\rightarrow \infty $.

Recall that in our case $h\left(x\right)=\sigma
\left(x\right)f_{S_0}\left(x\right)$ and the integral
\begin{align*}
&\int_{-\infty
}^{x } \sigma \left(y\right) f_{S_0}\left(y\right)\;{\rm
  d} \left[\frac{\zeta_T \left(S_0,y \right)}{2f_{S_0}\left(y\right)}\right]\\
&\qquad =\frac{1}{\sqrt{T}}\int_{0}^{T}\1_{\left\{X_t\leq x  \right\}}\;{\rm d}W_t
  -\frac{1}{\sqrt{T}}\int_{X_0}^{X_T}\frac{\1_{\left\{y\leq x\right\}}}{\sigma
    \left(y\right)}\;{\rm d}y.
\end{align*}

Therefore,
\begin{align*}
\hat \delta _T\left(X_T\right)&=\int_{-\infty }^{\infty }
\left(\frac{1}{\sqrt{T}}\int_{0}^{T}\1_{\left\{X_t\leq x
  \right\}}\;{\rm d}W_t -\frac{H\left(x,X_0,X_T\right)
}{\sqrt{T}}\right)^2{\rm d}F_{S_0}\left(x\right)\\
& =\int_{-\infty}^{\infty } \left(\frac{1}{\sqrt{T}}\int_{0}^{T}\1_{\left\{X_t\leq
  x \right\}}\;{\rm d}W_t \right)^2{\rm
  d}F_{S_0}\left(x\right)+O\left(\frac{1}{\sqrt{T}}\right),
\end{align*}
where we put
$$
H\left(x,y,z\right)=\int_{y}^{z}\frac{\1_{\left\{v\leq x\right\} }}{\sigma
    \left(v\right)}\;{\rm d}v.
$$
Note that the estimate $\Ex_{S_0}H\left(x,X_0,X_T\right)^2<C $ follows
directly from the assumption ${\cal A}_0$ (see \cite{Kut04} for details).

By  the law of large numbers
$$
\frac{1}{T}\int_{0}^{T}\1_{\left\{X_t\leq x
  \right\}}\;{\rm d}t\longrightarrow F_{S_0}\left(x\right).
$$
Hence, by the central limit theorem  the stochastic integral
\begin{align*}
I_T\left(x\right)&=\frac{1}{\sqrt{T}}\int_{0}^{T}\1_{\left\{X_t\leq x
  \right\}}\;{\rm d}W_t\\
&=\frac{1}{\sqrt{T}}\int_{0}^{T}\1_{\left\{ F_{S_0}\left(X_t\right)\leq F_{S_0}\left(x\right)
  \right\}}\;{\rm d}W_t\equiv J_T\left(F_{S_0}\left(x\right) \right)
\end{align*}
converges  
\begin{equation}
\label{fd}
\left(I_T\left(x_1\right) ,\ldots,I_T\left(x_k\right) \right)\Longrightarrow \left(I\left(x_1\right) ,\ldots,I\left(x_k\right)\right) ,
\end{equation}
where
$$
I\left(x\right)=w\left({F_{S_0}\left(x\right)}\right).
$$
Further, ($x_1<x_2$)
\begin{align}
\label{smo}
\Ex_{S_0}\left|I_T\left(x_2\right)- I_T\left(x_1\right) \right|^2=F_{S_0}\left(x_2\right)-F_{S_0}\left(x_1\right)\leq C\left|x_2-x_1\right|
\end{align}
 because under coondition ${\cal A}_0$ the density $f_{S_0}\left(\cdot \right)
 $  is a bounded function.

The properties \eqref{fd} and \eqref{smo} yield the convergence
$$
\int_{-\infty }^{\infty } \left[I_T\left(x\right)\right]^2\,{\rm
  d}F_{S_0}\left(x\right)=\int_{0}^{1 } \left[J_T\left(t\right)\right]^2\,{\rm
  d}t \Longrightarrow \int_{0}^{1 } w\left(t\right)^2\,{\rm  d}t
$$
(Theorem A22, \cite{IH81}).
 When we know the form of the statistic $\delta _T$ we can
construct another goodness of fit test with the same asymptotic properties   as follows.
Let us introduce the statistic
\begin{equation}
\label{main}
\delta _T^*=\int_{-\infty }^{\infty }\left( \frac{1}{\sqrt{T}}\int_{0}^{T}
  \frac{\1_{\left\{X_t<x\right\}}}{\sigma \left(X_t\right)}\left[{\rm
  d}X_t-S_0\left(X_t\right){\rm d}t\right] \right)^2{\rm d}F_{S_0}\left(x\right)
\end{equation}
and the constant $c_\varepsilon $:
$$
\Pb\left\{\int_{0}^{1 } w\left(t\right)^2\,{\rm  d}t>c_\varepsilon \right\}.
$$
We assume  that under the nonparametric alternative
$$
{\cal H}_1\quad \quad :\qquad \qquad \Ex_S\left(\frac{S\left(\xi
  \right)-S_0\left(\xi \right)}{\sigma \left(\xi \right)}\right)^2>0
$$
 the function
$S\left(\cdot \right)$ satisfies the conditions ${\cal ES}$ and ${\cal
  A}_0$. Here $\xi $ is the random variable with the density
 $f_S\left(x\right)$. Therefore the observed process is ergodic with the invariant
 density $f_S\left(\cdot \right)$.

Then for the test $\psi _T^*=\1_{\left\{\delta _T^*>c_\varepsilon  \right\}}$
we have the following result.
\begin{proposition}
\label{P1}
The test $\psi _T^*\in{\cal K}_\varepsilon $ and is consistent against any
fixed alternative ${\cal H}_1 $.
\end{proposition}
{\bf Proof.} Under hypothesis ${\cal H}_0$  we have
$$
\delta _T^*=\int_{-\infty }^{\infty }\left( \frac{1}{\sqrt{T}}\int_{0}^{T}
  {\1_{\left\{X_t<x\right\}}}\;{\rm
  d}W_t \right)^2{\rm d}F_{S_0}\left(x\right)
$$
and it converges to the following limit
$$
\delta _T^*\Longrightarrow \int_{0}^{1}w\left(t\right)^2{\rm d}t.
$$
Therefore, the test $\psi _T^*=\1_{\left\{\delta _T^*>c_\varepsilon \right\}}$
 belongs to ${\cal K}_\varepsilon $.

 The consisteny follows from  standard arguments as follows. Under alternative ${\cal H}_1$ we can write
\begin{align*}
&\frac{1}{\sqrt{T}}\int_{0}^{T}
  \frac{\1_{\left\{X_t<x\right\}}}{\sigma \left(X_t\right)}\left[{\rm
  d}X_t-S_0\left(X_t\right){\rm d}t\right]=\frac{1}{\sqrt{T}}\int_{0}^{T}
  {\1_{\left\{X_t<x\right\}}}{\rm
  d}W_t\\
&\qquad \quad +\frac{1}{\sqrt{T}}\int_{0}^{T}
  \frac{\1_{\left\{X_t<x\right\}}}{\sigma \left(X_t\right)}\left[S\left(X_t\right)-S_0\left(X_t\right)\right]{\rm d}t.
\end{align*}

The first (stochastic) integral is  asymptotically normal
$$
I_T\left(x\right)=\frac{1}{\sqrt{T}}\int_{0}^{T} {\1_{\left\{X_t<x\right\}}}{\rm
  d}W_t\Longrightarrow w\left(F_S\left(x\right)\right)\sim{\cal N}\left(0, F_S\left(x\right)\right)
$$
and for the second we have by the law of large numbers
$$
M_T\left(x\right)=\frac{1}{{T}}\int_{0}^{T}
  \frac{\1_{\left\{X_t<x\right\}}}{\sigma
    \left(X_t\right)}\left[S\left(X_t\right)-S_0\left(X_t\right)\right]{\rm
    d}t \longrightarrow  M\left(x\right).
$$
Here
$$
M\left(x\right)=\Ex_S\left(\frac{\1_{\left\{\xi <x\right\}}}{\sigma
    \left(\xi \right)}\left[S\left(\xi \right)-S_0\left(\xi \right)\right]\right).
$$

 Denote
$$
\left\|h \left(\cdot \right)\right\|_0^2=\int_{-\infty }^{\infty
}h\left(x\right)^2{\rm d}F_0\left(x\right) .
$$
Then we can write
\begin{align*}
\Pb_S\left( \delta _T^*>c_\varepsilon
\right)&=\Pb_S\left(\left\|I_T\left(\cdot \right)+\sqrt{T}
M_T\left(\cdot \right)\right\|_0>\sqrt{c_\varepsilon} \right)\\
&\geq \Pb_S\left(\sqrt{T}\left\|
M_T\left(\cdot \right)\right\|_0-\left\|I_T\left(\cdot
\right)\right\|_0>\sqrt{c_\varepsilon} \right).
\end{align*}
Therefore, if $\left\|M_T\left(\cdot \right)\right\|_0\rightarrow
\left\|M\left(\cdot \right)\right\|_0>0$ then the test is consistent.

The condition $\left\|M\left(\cdot \right)\right\|_0=0$ implies equality
$$
\int_{-\infty }^{x}\frac{\left[S\left(y\right)-S_0\left(y\right)\right]}{\sigma
    \left(y\right)}f_S\left(y\right){\rm d}y=0\quad {\rm for \;\;\; all}\quad
x\in R.
$$
Hence
$$
\frac{S\left(x\right)-S_0\left(x\right)}{\sigma
    \left(x\right)}f_S\left(x\right)=0\quad {\rm for \;\;almost \;\; all}\quad
x\in R,
$$
which is equivalent to the equality $S\left(x\right)=S_0\left(x\right)$ for almost  all
$x$. This contradicts the definition of the  alternative.

We do not consider here the ADF test based on empirical distribution function,
because the derivation of $\hat\eta _T\left(x\right)$  
reduces to the test  based on $\hat\zeta _T\left(x\right)$, for which we have already suggested a solution. 

{\bf Remark.} Note that the central statistic in \eqref{main} coincides with
the statistic used in the Kolmogorov-Smirnov type test studied \cite{NN}.  It
is interesting to note that another Kolmogorov-Smirnov test based on
empirical density $\zeta _T\left(x\right)=\sqrt{T}\left(\hat
f_T\left(x\right)-f_{S_0}\left(x\right)\right)$ was studied in
\cite{Kut04}. To compare these two tests we put $\sigma \left(x\right)\equiv
1$. Then the test in \cite{NN} is $\bar\psi_T\left(X^T\right)=\1_{\left\{
  \bar\delta _T>c_\varepsilon \right\} } $ with
\begin{equation}
\label{nn}
\bar\delta _T=\sup_x{\frac{1}{\sqrt{T}}
}\left|\int_{0}^{T}\1_{\left\{X_t<x\right\}} {\rm
    d}X_t-   \int_{0}^{T}\1_{\left\{X_t<x\right\}}  S_0\left(X_t\right){\rm
  d}t \right| .
\end{equation}
It is shown that  under hypothesis
$$
\bar\delta _T\Longrightarrow \sup_{0\leq t\leq 1}\left|w\left(t\right)\right|.
$$
Therefore this  test is ADF.

  The test proposed in  \cite{Kut04} is
$$
\delta _T^\circ=\sup_x\left|\zeta_T\left(x\right) \right|,\qquad
\psi_T^\circ\left(X^T\right)= \1_{\left\{\delta _T^\circ>c_\varepsilon  \right\}}.
$$

It is shown  that this statistic is asymptotically equivalent to the statistic
\begin{equation}
\label{yk}
\delta _T^+=\sup_x{\frac{2}{\sqrt{T}}
}\left|\int_{0}^{T}{\1_{\left\{X_t<x\right\}}}{\rm
    d}X_t-   \Ex_{S_0}\int_{0}^{T}{\1_{\left\{X_t<x\right\}}}   S_0\left(X_t\right){\rm d}t \right|
\end{equation}
and it converges to the limit
$$
\delta _T^+\Longrightarrow \sup_x\left|\Phi \left(x\right)\right|f_0\left(x\right).
$$
The comparison of \eqref{nn} and \eqref{yk} shows the advantage of \eqref{nn}
because it is ADF.

\end{document}